# Exponential Convergence through Linear Finite Element Discretization of Stratified Subdomains

Murthy N. Guddati [1*], Vladimir Druskin[2], Ali Vaziri Astaneh[1]

[1]*Department of Civil Engineering, North Carolina State University, Raleigh, NC 27695-7908, USA*

[2]*Schlumberger Doll Research, Cambridge, MA 02139, USA*

**Abstract**

Motivated by problems where the response is needed at select localized regions in a large computational domain, we devise a novel finite element discretization that results in exponential convergence at pre-selected points. The two key features of the discretization are (a) use of midpoint integration to evaluate the contribution matrices, and (b) an unconventional bending of the mesh into complex space. Named complex-length finite element method (CFEM), the technique is linked to Padé approximants that provide exponential convergence of the Dirichlet-to-Neumann maps and thus the solution at specified points in the domain. Exponential convergence facilitates drastic reduction in the number of elements. This, combined with sparse computation associated with linear finite elements, results in significant reduction in the computational cost. The paper presents the basic ideas of the method as well as illustration of its effectiveness for a variety of problems involving Laplace, Helmholtz and elastodynamic equations.

## 1. Introduction

Conventional domain-based methods such as finite element and finite difference techniques obtain the solution over the entire domain. While such approach is appropriate in many problems, there are several situations where the response is needed only in a few small regions of interest. Some example include: (a) reservoir modeling where the response at injection and production wells are of utmost importance, (b) forward modeling in the context of nondestructive testing and system identification, where the goal is to match the response of the system at sparse discrete points in the domain and (c) structural acoustics, where the acoustic signature is not needed at all the points, but at distinct locations in the far-field. Most of these problems involve significant computational expense and it would be desirable to reduce the computational cost, if it can come at the expense of not computing the response in the rest of the domain. With this motivation, this paper presents an unconventional finite element method that provides high accuracy at prescribed points in the domain. In this paper, we treat the special but important class of problems involving large regular (e.g. layered) subdomains, where the actual solution inside these subdomains may not be of interest, but rather the effect of these subdomains on the solution in the remainder of domain. Specifically, we show that a special grid with midpoint-integrated linear finite elements can result in exponential convergence of the solution on the edges of layered subdomains.

* Corresponding author. Tel.: +1 919 515 7699; fax: +1 919 515 7908
  E-mail address: murthy.guddati@ncsu.edu

Exponential convergence is typically achieved with the help of spectral methods where the field variable is discretized with Fourier basis [1, 2]. Unfortunately, spectral methods in general render the computation global. On the other hand, regular finite element and finite difference methods involve more efficient sparse computation, but the convergence is only algebraic. It was discovered that exponential convergence can be obtained with sparse computation, provided that the accuracy is needed only at specific points in the domain [3, 4]. By linking finite-difference approximation to rational approximation of the Dirichlet-to-Neumann (DtN) map, exponential convergence is achieved at the edges of sub-domains discretized with specially devised finite difference grids. The basic idea is to obtain optimal rational approximation of the *one-sided* DtN map (with Dirichlet or Neumann condition applied on the other edge), and translate the approximation to an equivalent finite difference grid. Since the grids are obtained through exponentially convergent optimal approximations of the DtN map, they are called *optimal finite difference grids* and result in exponential convergence at the edges of the sub-domains. The main limitation of this method is that two distinct finite difference grids are needed, one when Dirichlet condition is applied on the other edge, and another for Neumann condition, indicating that the grids cannot be used directly for *two-sided* problems, which would be the building block for multi-domain problems. Optimal grids can be applied to two-sided problems *indirectly*, through splitting the solution into odd and even parts, devising two separate grids for each part, and using them in a completely overlapping function. This idea is extended to multi-dimensions, but the computation becomes rather cumbersome, requiring increasing number of overlapping grids (four for two-dimensional problems and eight for three-dimensions) [4-6].

In this paper, we introduce a simpler alternative to optimal finite difference grids and eliminate the need for overlapping grids. The key to the proposed method is the observation that *midpoint integrated linear finite elements* facilitate fixed-point iteration for the DtN map of the half-space, i.e. adding these elements to a half-space does not alter the DtN map of the composite half-space [7-9]. We show that this property eliminates the need for multiple overlapping grids and provides *exponential convergence of the DtN map for the two-sided problem with a single grid*. This makes the implementation attractive and the computation can be performed by a simple modification of existing finite element codes. The only complication is that the finite element mesh needs to be bent out of the real space, making the element lengths complex-valued. This feature necessitates complex arithmetic and could contribute to an increase in the computational cost. However, this increase is not an issue as the proposed method needs very coarse grid, with number of elements typically less by an order of magnitude than that for regular finite element discretization. Moreover, in many cases including time-harmonic wave propagation and electromagnetism, the original computation involves complex arithmetic and the bending of the mesh into complex space does not add any further complications.

This paper focuses on the derivation of the method along with numerical examples to illustrate its efficiency. The treatment is limited to the basic development of the method, with further extensions left for future research. The outline of the development is as follows. We start with the two-dimensional model problem given in Fig. 1 and show that it can be reduced to a set of one-dimensional two-point boundary value problems through semi-discretization in the vertical direction (Section 2). We then try to obtain a grid that simultaneously approximates the DtN map for all the one-dimensional problems (Section 3). This is achieved by reformulating midpoint integrated linear finite elements, described in Section 3.1, as Crank-Nicolson discretization of an equivalent first-order system (Section 3.2). Exponential convergence is achieved by choosing the parameters by linking the Crank-Nicolson discretization with Padé approximants (Section 3.3). The result is a finite element mesh with complex element lengths, which is discussed in Section 3.3. In Section 4 we show the applicability of the present method application to general vector (elastic) equations. Section 5 contains numerical examples illustrating the exponential convergence and practical use of the proposed method. The paper is concluded in Section 6 with closing remarks.

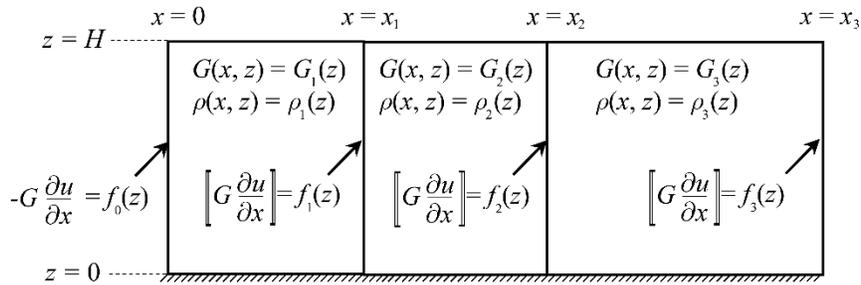

**Fig. 1.** A schematic of the model problem. An elastic layer is deforming out of plane under external excitation $f_{0,1,2,3}$ at the vertical interfaces. The segments bounded by these interfaces are stratified in the vertical direction. The goal is to obtain the responses at the interfaces.

## 2. Preliminaries

*2.1. Problem statement*

For the sake of focused discussion, consider the (artificial) model problem in Fig. 1, where an elastic layer is made up of three segments, each segment individually stratified in the vertical direction. Excitation is applied only at the vertical edges and interfaces and the goal is to obtain the response only at these locations. We consider the anti-plane shear deformation governed by the Helmholtz equation (Laplace equation being a special case when frequency $\omega = 0$):

$$-\frac{\partial}{\partial z}\left(G(x,z)\frac{\partial w}{\partial z}\right) - \frac{\partial}{\partial x}\left(G(x,z)\frac{\partial w}{\partial x}\right) - \rho(x,z)\omega^2 w = 0, \qquad (1)$$

where $w$ is the anti-plane deformation, $G$ is the shear modulus and $\rho$ is the density. Since the forces are applied only at the vertical interfaces and the response is needed only at these locations, it is sufficient to solve the coupled problem obtained by assembling the Dirichlet-to-Neumann (DtN) maps of individual segments. We thus focus on obtaining an accurate estimate of the DtN map for a single segment, which reduces the problem to:

$$\textbf{Problem 1:} \begin{cases} \text{Obtain the relationship DtN: } w\big|_{x=0,L} \to \dfrac{dw}{dx}\bigg|_{x=0,L}, \text{ with } w \text{ satisfying:} \\[6pt] -\dfrac{\partial}{\partial z}\left(G(z)\dfrac{\partial w}{\partial z}\right) - \dfrac{\partial}{\partial x}\left(G(z)\dfrac{\partial w}{\partial x}\right) - \rho(z)\omega^2 w = 0, \quad \text{for } \begin{array}{l} 0 < x < L, \\ 0 < z < H, \end{array} \\[6pt] w\big|_{z=0} = 0, \qquad \dfrac{\partial w}{\partial z}\bigg|_{z=H} = 0. \end{cases} \qquad (2)$$

*2.2. Dimensional reduction*

Our ultimate goal is to obtain a finite element mesh that can yield accurate results at the vertical interfaces. Considering that accuracy is needed along the entire length of every interface, we propose to perform $z$ discretization in a conventional manner and independent of $x$ coordinate. It is also natural to make $x$ discretization independent of $z$, resulting in a two-dimensional tensor-product mesh. For the purposes of analyzing the method, the tensor product discretization can be viewed as discretization in $z$ followed by discretization in $x$. In this subsection, we focus on the discretization in $z$ and the resulting dimensionally reduced differential equation in $x$.

We start with obtaining the weak form in $z$ by (a) multiplying the governing equation by virtual displacement by $\delta w$, (b) performing integration by parts only in $z$, and (c) applying the boundary

conditions at $z = 0, H$. This results in:

$$\int_{z=0}^{z=L} \left( \frac{\partial \delta w}{\partial z} G(z) \frac{\partial w}{\partial z} + \delta w G(z) \frac{\partial^2 w}{\partial x^2} - \delta w \omega^2 \rho(z) w \right) dz = 0 \quad \text{for all } \delta w, \quad 0 < x < L. \tag{3}$$

We apply Bubnov-Galerkin method, i.e. use the following approximations,

$$w \approx \mathbf{N}_z(z)\mathbf{W}(x) \quad \text{and} \quad \delta w \approx \mathbf{N}_z(z)\delta\mathbf{W}(x), \tag{4}$$

where $\mathbf{N}_z$ is the matrix of interpolation functions (shape function matrix) and $\mathbf{W}$ and $\delta\mathbf{W}$ are vectors of discretized field and virtual field variables respectively. We emphasize that $\mathbf{W}$ and $\delta\mathbf{W}$ are functions of $x$ as the discretization is performed only in the $z$ direction; the goal is to solve for unknown $\mathbf{W}(x)$. After substituting (4) in (3) and performing appropriate manipulations, we obtain the discrete form of the governing equation:

$$\mathbf{RW} - \mathbf{H}\frac{d^2\mathbf{W}}{dx^2} - \omega^2\mathbf{MW} = \mathbf{0}, \quad \text{for } 0 < x < L, \tag{5}$$

where

$$\mathbf{R} = \int_{z=0}^{z=H} \frac{\partial \mathbf{N}_z}{\partial z}^T G \frac{\partial \mathbf{N}_z}{\partial z} dz, \quad \mathbf{H} = \int_{z=0}^{z=H} \mathbf{N}_z^T G \mathbf{N}_z dz, \quad \mathbf{M} = \int_{z=0}^{z=H} \mathbf{N}_z^T \rho \mathbf{N}_z dz. \tag{6}$$

Eq. (5) can be easily decoupled into a set of differential equations,

$$-\frac{\partial^2 u_l}{\partial x^2} + \lambda_l u_l = 0 \quad \text{for } 0 < x < L, \quad \lambda_{\min} \leq \lambda_l \leq \lambda_{\max}, \tag{7}$$

where $\lambda_l$ are the eigenvalues of $\mathbf{H}$ with respect to $\mathbf{R} - \omega^2\mathbf{M}$. $u_l$ are the weights of the corresponding eigenvectors (modal participation factors). Note that $\lambda_l$ are real and bounded because all the matrices are symmetric and $\mathbf{A}$ is positive definite.

With the above decoupling procedure, the original problem of obtaining the two-dimensional DtN map simplifies to finding the DtN maps for a set of one-dimensional boundary value problems:

$$\textbf{Problem 2: } \begin{cases} \text{Obtain the relationship DtN: } u|_{x=0,L} \to \left.\dfrac{du}{dx}\right|_{x=0,L}, \text{ with } u \text{ satisfying:} \\ -\dfrac{\partial^2 u}{\partial x^2} + \lambda u = 0 \quad \text{for } 0 < x < L, \quad \lambda_{\min} \leq \lambda \leq \lambda_{\max}. \end{cases} \tag{8}$$

Note that the subscript $l$ is removed for the sake of simplicity in presentation. For the special case of homogeneous media, the shape functions $\mathbf{N}_z$ can be chosen according to Fourier expansion, resulting in the classical spectral element method. This would automatically result in decoupled equations in (7) without the need for eigendecomposition. Nevertheless, we use general discretization in the $z$ direction, to not limit the applicability to homogeneous media, or to a scalar equation. Note that eigen decomposition in (7) is performed only for the purposes of analyzing the method; it is never performed in the actual computation.

## 2.3. A closer look at the DtN map

The desired DtN map is a $2 \times 2$ matrix,

$$\begin{Bmatrix} -v_0 \\ v_L \end{Bmatrix} = \underbrace{\begin{bmatrix} K_{00} & K_{0L} \\ K_{L0} & K_{LL} \end{bmatrix}}_{=\mathbf{K}} \begin{Bmatrix} u_0 \\ u_L \end{Bmatrix}, \tag{9}$$

where $v = \partial u / \partial x$ and the subscripts represent the locations at which the quantities are evaluated. The negative sign for $v_0$ is used since the outside normal at $x = 0$ is in the negative $x$ direction. Considering the geometric symmetry of the domain as well as the symmetry of the operator, the DtN map can be written as,

$$\mathbf{K} = \begin{bmatrix} K_{diag} & K_{off} \\ K_{off} & K_{diag} \end{bmatrix}. \tag{10}$$

For a fixed $\lambda$, the exact DtN map can be easily obtained by solving the two-point boundary value problem in (8), and is given by,

$$\mathbf{K}_{exact} = \begin{bmatrix} K_{diag} & K_{off} \\ K_{off} & K_{diag} \end{bmatrix} = \frac{\sqrt{\lambda}}{\sinh(\sqrt{\lambda}L)} \begin{bmatrix} \cosh(\sqrt{\lambda}L) & -1 \\ -1 & \cosh(\sqrt{\lambda}L) \end{bmatrix}. \tag{11}$$

Note that the exact DtN map satisfies the *fixed-point property*, which is defined as follows: when $\mathbf{K}_{exact}$ is augmented with the exact DtN map of a half-space $\left(\mathbf{K}_{halfspace} = \sqrt{\lambda}\right)$ and the interior node is eliminated, exact half-space stiffness is preserved at the exterior node (this follows from the simple physical argument that when a layer is added to a half-space with same material properties, we obtain the same half-space). In other words, fixed point property states that, when a layer $(x_0, x_1)$ is augmented with a half-space $(x_1, \infty)$, the Neumann data required at $x_0$ to generate a Dirichlet data of $u_0$ is $K_{halfspace} u_0$, i.e.,

$$\begin{bmatrix} K_{diag} & K_{off} \\ K_{off} & K_{diag} + K_{halfspace} \end{bmatrix} \begin{Bmatrix} u_0 \\ u_1 \end{Bmatrix} = \begin{Bmatrix} K_{halfspace} u_0 \\ 0 \end{Bmatrix}. \tag{12}$$

It is easy to show that the above definition is equivalent to,

$$K_{halfspace} = K_{diag} - \frac{K_{off}^2}{\left(K_{diag} + K_{halfspace}\right)}. \tag{13}$$

Noting that $K_{halfspace} = \sqrt{\lambda}$, it is easy to see that fixed-point property takes a simple form,

$$K_{diag}^2 - K_{off}^2 = \lambda. \tag{14}$$

While it is obvious that the exact DtN map satisfies the fixed-point property, approximate DtN maps do not necessarily satisfy this property. The discretization proposed here, however, results in a DtN map that does satisfy the fixed point property. This is the key to the current development, as discussed in Section 3.1.

The exact DtN map in (11), while useful for any given $\lambda$, cannot be used for solving the original 2D problem. The 2D problem requires a grid that would *simultaneously* work for *all* $\lambda$. One option is to use the standard equidistant finite element or finite difference grid; such as grid results in algebraic convergence of the solution over the entire domain, and thus the DtN map. Our goal involves accurate approximation of only the DtN map. It may be possible to obtain more accuracy for the DtN map, possibly at the expense of the solution accuracy in the interior. Druskin and coworkers [3, 4, 10-12] have successfully utilized this idea to obtain exponential convergence of one-sided DtN map, i.e., the DtN map at one end of the domain with the boundary condition specified at the other end. Utilizing rational approximation theory of Stieltjes functions and its link to staggered-grid finite difference method, they obtained two independent grids that provide exponential convergence of the one-sided DtN map, one for Dirichlet boundary condition at the other end, and the other for Neumann boundary condition. They also extended the method to the two-sided problem, but the extension is not straightforward. They split the solution into odd and even parts and use two separate grids to approximate the corresponding DtN maps. Simultaneous solution of two-point boundary value problem requires the use of these grids in a completely overlapping fashion, making the implementation cumbersome, especially when extended to higher dimensions (we would need to use 4 overlapping grids for 2D and eight for 3-D) [4-6]. Therefore, it is desirable to obtain a single mesh that would simultaneously approximate all the elements of the $2\times 2$ DtN map with equal accuracy. As will be shown in the rest of the paper, the midpoint-integrated linear finite elements proposed by Guddati and coworkers [7-9] leads to such approximation. The key property of these elements is that they preserve the fixed-point property of the exact DtN map, as defined in the previous section. We utilize this property and develop linear finite element discretization that provides simultaneous exponential convergence for all the elements of the DtN map. The resulting finite element mesh is *not* a standard discretization of the real domain $(0, L)$, but is bent into the complex space (the element lengths are complex-valued). We term this the *complex-length finite element method* (*CFEM*). The formulation of the method is given in the following section.

## 3. The complex-length finite element method (CFEM)

We start this section with a brief overview of mid-point integrated linear finite elements and their special properties, followed by a detailed analysis, eventually resulting in an exponentially convergent finite element mesh.

## 3.1. Linear finite elements with midpoint integration

Consider that the domain $(0, L)$ is discretized into non-overlapping finite elements of lengths $L_j$, $j = 1, \cdots, n$, with $\sum L_j = L$, with nodes located at $x_{0,1,\cdots,n}$. Restricting the discussion to $j$ th element, the weak form becomes,

$$\int_{x_{j-1}}^{x_j} \left( \frac{\partial \delta u}{\partial x} \frac{\partial u}{\partial x} + \delta u \lambda u \right) dx = \delta u_j v_j - \delta u_{j-1} v_{j-1}, \quad (15)$$

where $\delta u$ is the variation of $u$. Utilizing Bubnov-Galerkin method with linear interpolation within the element and after performing the transformation $\bar{x} = x - x_{j-1}$, we obtain,

$$\begin{Bmatrix} \delta u_{j-1} \\ \delta u_j \end{Bmatrix}^T \left[ \int_0^{L_j} \left( \begin{Bmatrix} -1/L_j \\ 1/L_j \end{Bmatrix} \begin{Bmatrix} -1/L_j \\ 1/L_j \end{Bmatrix}^T + \begin{Bmatrix} 1-\bar{x}/L_j \\ \bar{x}/L_j \end{Bmatrix} \lambda \begin{Bmatrix} 1-\bar{x}/L_j \\ \bar{x}/L_j \end{Bmatrix}^T \right) d\bar{x} \right] \begin{Bmatrix} u_{j-1} \\ u_j \end{Bmatrix} = \begin{Bmatrix} \delta u_{j-1} \\ \delta u_j \end{Bmatrix}^T \begin{Bmatrix} -v_{j-1} \\ v_j \end{Bmatrix}. \quad (16)$$

Next step is to evaluate the integral using the midpoint rule. Performing such integration and eliminating the virtual displacements, we obtain the DtN map (stiffness relation) for a single element,

$$\begin{Bmatrix} -v_{j-1} \\ v_j \end{Bmatrix} = \begin{bmatrix} \left( \frac{1}{L_j} + \frac{\lambda L_j}{4} \right) & \left( -\frac{1}{L_j} + \frac{\lambda L_j}{4} \right) \\ \left( -\frac{1}{L_j} + \frac{\lambda L_j}{4} \right) & \left( \frac{1}{L_j} + \frac{\lambda L_j}{4} \right) \end{bmatrix} \begin{Bmatrix} u_{j-1} \\ u_j \end{Bmatrix}. \quad (17)$$

In the above DtN map, $K_{\text{diag}}^2 - K_{\text{off}}^2 = \lambda$, which indicates the fixed-point property (see (14)). Since all the finite elements in the mesh satisfy this property, the two-point DtN map of the entire mesh would also satisfy the fixed-point property. An important consequence is that, a mesh that results in exponential convergence of the *diagonal* element of the DtN map, also results in exponential convergence of the *off-diagonal* element (this is because, $K_{\text{diag, approx}}^2 - K_{\text{off, approx}}^2 = \lambda = K_{\text{diag, exact}}^2 - K_{\text{off, exact}}^2$), and thus exponential convergence of the entire two-sided DtN map.

To obtain the DtN map of the discretized domain, we could assemble the $2 \times 2$ DtN maps of the individual elements and eliminate the interior degrees of freedom in the sense of taking Schur complement. However, this process renders the computation global, making it cumbersome to analyze. To facilitate easier analysis, we turn to an equivalent propagator matrix approach, by converting the boundary value problem into an initial value problem.

## 3.2. An alternative view: first-order form

With simple substitution, we can verify that the DtN relation in (17) is equivalent to,

$$\begin{Bmatrix} \dfrac{u_j - u_{j-1}}{L_j} \\ \dfrac{\bar{v}_j - \bar{v}_{j-1}}{L_j} \end{Bmatrix} = \begin{bmatrix} 0 & \sqrt{\lambda} \\ \sqrt{\lambda} & 0 \end{bmatrix} \begin{Bmatrix} \dfrac{u_j + u_{j-1}}{2} \\ \dfrac{\bar{v}_j + \bar{v}_{j-1}}{2} \end{Bmatrix}, \tag{18}$$

where $\bar{v} = v/\sqrt{\lambda}$. It can be immediately seen that the above equation is the Crank-Nicolson discretization of first order form of differential equation (7), i.e.,

$$\frac{\partial}{\partial x}\begin{Bmatrix} u \\ \bar{v} \end{Bmatrix} = \begin{bmatrix} 0 & \sqrt{\lambda} \\ \sqrt{\lambda} & 0 \end{bmatrix}\begin{Bmatrix} u \\ \bar{v} \end{Bmatrix} \quad \text{for } 0 < x < L, \quad \lambda_{min} \leq \lambda \leq \lambda_{max}. \tag{19}$$

Eq. (18) can also be written in the form,

$$\begin{Bmatrix} u_j \\ \bar{v}_j \end{Bmatrix} = \mathbf{P}_j \begin{Bmatrix} u_{j-1} \\ \bar{v}_{j-1} \end{Bmatrix}, \tag{20}$$

where $\mathbf{P}_j$ is the *propagator matrix* associated with the finite element, which is inherently connected to the DtN map in (17). The advantage of the propagator matrix approach is that the propagator matrix of the complete interval $(0, L)$, $\mathbf{P}$, is obtained by simply multiplying the propagator matrices of individual intervals:

$$\begin{Bmatrix} u_L \\ \bar{v}_L \end{Bmatrix} = \mathbf{P}_n \begin{Bmatrix} u_{n-1} \\ \bar{v}_{n-1} \end{Bmatrix} = \mathbf{P}_n \mathbf{P}_{n-1} \begin{Bmatrix} u_{n-2} \\ \bar{v}_{n-2} \end{Bmatrix} = \cdots = \underbrace{\mathbf{P}_n \mathbf{P}_{n-1} \cdots \mathbf{P}_2 \mathbf{P}_1}_{\mathbf{P}} \begin{Bmatrix} u_0 \\ \bar{v}_0 \end{Bmatrix}. \tag{21}$$

The matrix $\mathbf{P}$ can also be written in the form of DtN relation in (9), implying that successful approximation of the propagator matrix would automatically result in successful approximation of the DtN map. Thus, the problem is reduced to *finding the mesh parameters $L_j$ that would result in exponential convergence of the propagator matrix $\mathbf{P}$ for various values of $\lambda$ on the real line.*

The problem of approximating the propagator matrix can be further simplified by decoupling the system of equations in (19). Through eigendecomposition, we have,

$$\frac{\partial}{\partial x}\begin{Bmatrix} \psi_1 \\ \psi_2 \end{Bmatrix} = \begin{bmatrix} +\sqrt{\lambda} & 0 \\ 0 & -\sqrt{\lambda} \end{bmatrix}\begin{Bmatrix} \psi_1 \\ \psi_2 \end{Bmatrix} \quad \text{for } 0 < x < L, \quad \lambda_{min} \leq \lambda \leq \lambda_{max}, \tag{22}$$

where the diagonal matrix contains the eigenvalues of the operator in (19), and $\psi_1, \psi_2$ are the weights associated with the corresponding eigenvectors. Thus the DtN approximation problem reduces to:

$$\textbf{Problem 3:} \begin{cases} \text{Find the mesh parameters } L_j \text{ that would result in exponentially} \\ \text{convergent propagator matrix for the interval } (0,L) \text{ for the equation:} \\ \dfrac{\partial \psi}{\partial x} = k\psi \quad \text{for all } k \in \left(-\sqrt{|\lambda_{\max}|}, \sqrt{|\lambda_{\max}|}\right) \cup \left(-i\sqrt{|\lambda_{\min}|}, i\sqrt{|\lambda_{\min}|}\right). \end{cases} \quad (23)$$

In the above, $k = \sqrt{\lambda}$ and $\lambda_{\min}$ is implicitly assumed to be negative, indicating that $k$ could be imaginary. It is instructive to note the physical meaning of the eigenvectors associated with the decoupled form in (22). The eigenvector associated with $\psi_1$ is $\{u \ \bar{v}\} \propto \{1 \ 1\}$, or equivalently, $\bar{v} = u$, which corresponds to an exponentially growing solution, $u = \psi_1 e^{\sqrt{\lambda}x}$. The other eigenvector is $\{u \ \bar{v}\} \propto \{1 \ -1\}$, which corresponds to a decaying solution, $u = \psi_2 e^{-\sqrt{\lambda}x}$. These are the solutions of the original second order differential equation (7). This point of view is important as our method is based on approximation of the exponential functions (see Section 3.3).

The fixed-point property of the midpoint-integrated linear finite elements can be proven more elegantly as follows. The impedance relation for the half-space $(0,\infty)$ is: $\partial u/\partial x|_{x=0} = \sqrt{\lambda} u_0$, or $\bar{v}_0 = u_0$ (this corresponds to the decaying solution). The fixed-point property (preservation of the half-space impedance relation) can be written as: $\bar{v} = u$ at $x = 0$ implies $\bar{v} = u$ at any $x$. The fixed-point property of the exact propagator can be illustrated as follows: Eq. (19) is invariant to the swapping of $u$ with $\bar{v}$; therefore, if the initial condition at $x = 0$ is swap-invariant ($\bar{v} = u$), the "final" condition at any $x$ must also be swap-invariant ($\bar{v} = u$). Fortunately, the swap-invariance of (19) is preserved through the discretization in (18), implying that the Crank-Nicolson propagator, and thus the midpoint-integrated finite element, satisfy the fixed-point property.

*3.3. Selection of the mesh parameters*

Considering that the exact solution to (23) is of the form $Ae^{kx}$, the propagator associated with the interval $(0,L)$ is simply,

$$P_{exact} = \frac{\psi_L}{\psi_0} = e^{kL}. \qquad (24)$$

The approximate propagator for the $j$ th interval is given by the Crank-Nicolson method,

$$\frac{\psi_j - \psi_{j-1}}{L_j} = k\left(\frac{\psi_j + \psi_{j-1}}{2}\right) \quad \Rightarrow \quad \psi_j = \underbrace{\left(\frac{1 + kL_j/2}{1 - kL_j/2}\right)}_{=P_j} \psi_{j-1}. \qquad (25)$$

The propagator for the entire interval $(0,L)$ is the product of all the propagators, i.e.,

$$P_{approx} = \prod_{j=1}^{n} P_j = \prod_{j=1}^{n} \left( \frac{1+kL_j/2}{1-kL_j/2} \right). \tag{26}$$

Note that $P_{approx}$ is a so-called relative approximant of $e^{kL}$, which is a rational function of the form $Q(k)/Q(-k)$, where $Q$ is a polynomial of degree $n$ with roots $2/L_j$. Relative approximants are very well-studied subject of rational approximation theory [13], and there is a variety of approximants convergent on the real axis at least exponentially with respect to $n$. Obviously, $Q(k)/Q(-k)$ must not have poles and residues on the real axis in order to be a good approximant of the exponent there; this is the reason why $L_j$ must be complex, implying that the mesh is bent into the complex space. Hence the current method is called the complex-length finite element method.

In this paper, we consider the relative (diagonal) Padé approximant matching the first $2n$ terms of Taylor expansion at $k=0$. In other words, the 0 to $2n^{th}$ derivatives of the approximant must match with that of the exact solution at $k=0$, i.e.,

$$\left. \frac{d^j P_{pade}}{dk^j} \right|_{k=0} = \left. \frac{d^j \left(e^{kL}\right)}{dk^j} \right|_{k=0}, \qquad j=0,\cdots,2n. \tag{27}$$

The reason for the above choice of Padé approximant is its simplicity. Even though the convergence is faster than exponential for any fixed real $k$, the convergence rate quickly deteriorates for large absolute values of $k$, which can create problems in approximating non-smooth PDE solutions. There are a number of approximants of the form $Q(k)/Q(-k)$ with better convergence properties for large $|k|$; we plan to consider such approximations in our future research.

For the specific Padé approximant considered in (27), the roots $2/L_j$ and thus the element lengths $L_j$ can be computed with a standard algorithms (see e.g. [14]). The resulting mesh has the following properties:

1. $L_j$ are independent of $k$ (or $\lambda$), indicating that the same mesh can be used for the range of complex $\lambda$ given in (8), implying that the mesh is applicable to the original 2D problem in (2).

2. $L_j$ come in complex conjugate pairs. It follows from the fact that the Padé approximant is real for any real $k$.

3. $\sum L_j = L$ (follows directly from evaluating (27) for $j=1$). This indicates that the end points match with the physical edges of the domain. This is an important property as the procedure for the sub-domain does not affect the geometry of the rest of the analysis domain.

4. A peculiar property of midpoint-integrated finite element grid is that the DtN map is invariant of the ordering of the elements. This implies that the mesh is not unique; we order the mesh so that it is symmetric and smooth (as smooth as possible) in the complex plane. Such a mesh is obtained by ordering the elements with increasing or decreasing phase angle of their lengths.

5. The mesh scales with the length of the domain and we can tabulate the element sizes relative to the domain size, i.e. $L_j/L$, for any given number of elements; this is done in Table 1 up to 16 elements. In general one can obtain $L_j/L = 2/x_j$ where $x_j$ are the roots of the $n^{th}$-degree polynomial given below (see e.g. [14]):

$$\sum_{j=0}^{n} \frac{(2n-j)!}{j!(n-j)!}(-x)^j = 0. \tag{28}$$

Also, following the order suggested in item 4 above, the meshes can also be graphically depicted by drawing the mesh in the complex plane (see Fig. 2). Separate figures are used for odd and even number of elements for the sake of clarity in presentation. With refinement, the mesh appears to converge to a specific curve in the complex space. We suspect that there is an analytical form for the converged mesh, linked to the asymptotic behavior of the poles of Padé approximants, similar to the one known for optimal approximants on $[0, \infty)$ [15].

6. It is known that the real parts of the roots $(2/L_j)$ are positive for the diagonal Padé approximants of exponential function [14, 16], implying that the real parts of $L_j$ are positive and the mesh is not bent outside the domain (or have any knots), which is another physically appealing aspect.

*3.4. Summary: the algorithm and properties*

The analysis of the model problem described in Fig. 1 can be solved using the following strategy:

1. Divide the problem into multiple sub-domains bounded by points where (a) the material properties are discontinuous, (b) loads are applied, or (c) the solution is desired.
2. Discretize each sub-domain using the complex-length finite element mesh obtained using the above procedure.
3. Combine all the sub-domain grids and perform finite element analysis on the entire domain. Utilize midpoint integration to evaluate the contribution matrices.
4. Extract the results from the sub-domain interfaces. The field variables in the interior points of the sub-domains have no immediate physical meaning.

**Table 1** Element lengths for varying number of elements ($n$). The domain size is assumed to be unity. For other cases, the element lengths simply scale with domain size.

| | | | |
|---|---|---|---|
| $n=1$ | 1.00000000000000 | $n=2$ | $0.50000000000000 \pm 0.28867513459481\,i$ |
| $n=3$ | $0.28468557688388 \pm 0.27159985141630\,i$<br>0.43062884623222 | $n=4$ | $0.18313248053143 \pm 0.23132522602625\,i$<br>$0.31686751946856 \pm 0.09488202514221\,i$ |
| $n=5$ | $0.12803667831541 \pm 0.19668213834621\,i$<br>$0.23485450871940 \pm 0.12209940763707\,i$<br>0.27421762593037 | $n=6$ | $0.09489061789607 \pm 0.16944514819433\,i$<br>$0.17914640739749 \pm 0.12594324946340\,i$<br>$0.22596297470643 \pm 0.04614135671779\,i$ |
| $n=7$ | $0.07338559568636 \pm 0.14811940741461\,i$<br>$0.14065739395847 \pm 0.12154781833235\,i$<br>$0.18538954553266 \pm 0.06776497788782\,i$<br>0.20113492964499 | $n=8$ | $0.05861791492234 \pm 0.13119236974564\,i$<br>$0.11325833004971 \pm 0.11445496413908\,i$<br>$0.15337794771885 \pm 0.07709430353886\,i$<br>$0.17474580730908 \pm 0.02713226173787\,i$ |
| $n=9$ | $0.04802049907890 \pm 0.11752570488080\,i$<br>$0.09316287173966 \pm 0.10679717990370\,i$<br>$0.12840441052931 \pm 0.08014721079992\,i$<br>$0.15100957026047 \pm 0.04281546509628\,i$<br>0.158805296783284 | $n=10$ | $0.04014472910062 \pm 0.10630697796687\,i$<br>$0.07802273616547 \pm 0.09940003581225\,i$<br>$0.10881874816902 \pm 0.07996015060428\,i$<br>$0.13078256453612 \pm 0.05163001172938\,i$<br>$0.14223122202876 \pm 0.01782841382104\,i$ |
| $n=11$ | $0.03412261657800 \pm 0.09695789626293\,i$<br>$0.06634486381072 \pm 0.09256005182226\,i$<br>$0.09329088025857 \pm 0.07811366645707\,i$<br>$0.11386072467713 \pm 0.05628616844255\,i$<br>$0.12678425930221 \pm 0.02942684799389\,i$<br>0.131193310746699 | $n=12$ | $0.02940803944815 \pm 0.08906181395662\,i$<br>$0.05715192456673 \pm 0.08635352530097\,i$<br>$0.08082582076929 \pm 0.07545289071966\,i$<br>$0.09976290741316 \pm 0.05839885662872\,i$<br>$0.11301395814524 \pm 0.03687019624343\,i$<br>$0.11983734965741 \pm 0.01259866487095\,i$ |
| $n=13$ | $0.02564318775369 \pm 0.08231355321087\,i$<br>$0.04978573946440 \pm 0.08076582096948\,i$<br>$0.07069390925651 \pm 0.07243833050796\,i$<br>$0.08799247146698 \pm 0.05894563008011\,i$<br>$0.10097469339377 \pm 0.04152986688324\,i$<br>$0.10902977699627 \pm 0.02144314210934\,i$<br>0.11176044333671 | $n=14$ | $0.02258550311646 \pm 0.07648569373967\,i$<br>$0.04379127631258 \pm 0.07574740948845\,i$<br>$0.06236027779351 \pm 0.06932300610809\,i$<br>$0.07811546938314 \pm 0.05852853926583\,i$<br>$0.09053178744825 \pm 0.04431195125450\,i$<br>$0.09911464460927 \pm 0.02760255210502\,i$<br>$0.10350104133675 \pm 0.00937153174338\,i$ |
| $n=15$ | $0.02006570730347 \pm 0.07140591920396\,i$<br>$0.03884638966705 \pm 0.07123849065223\,i$<br>$0.05542991160544 \pm 0.06624511011050\,i$<br>$0.06977487505357 \pm 0.05752432367232\,i$<br>$0.08149231310591 \pm 0.04581962272804\,i$<br>$0.09018490369015 \pm 0.03183750267162\,i$<br>$0.09553514496841 \pm 0.01630985840134\,i$<br>0.09734150921194 | $n=16$ | $0.01796266496341 \pm 0.06694162366630\,i$<br>$0.03471809507502 \pm 0.06717964470869\,i$<br>$0.04960794930403 \pm 0.06327803486714\,i$<br>$0.06268398209122 \pm 0.05617188543383\,i$<br>$0.07365957371625 \pm 0.04645858465651\,i$<br>$0.08221667350938 \pm 0.03468717833984\,i$<br>$0.08808374597041 \pm 0.02141958057287\,i$<br>$0.09106731537025 \pm 0.00724175169196\,i$ |

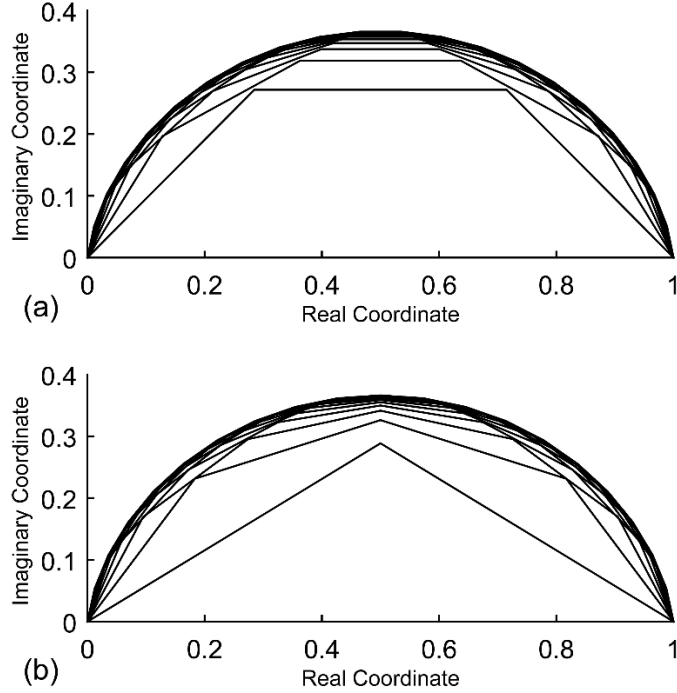

**Fig. 2.** Finite element meshes for (a) odd number of elements and (b) even number of elements.

It is also worthwhile mentioning that the discussed method of CFEM satisfies the following important properties:

1. Our method can be considered as a standard finite element method with complex coordinate stretching. Thus, it shares similarities with perfectly matched layers (PML) used in unbounded domain modeling, which can be viewed in terms of complex coordinate stretching [17]. The difference is that PML stretching leaves invariant one boundary point, whereas both boundary points are invariant in our method. PML finite elements lengths do not come in conjugate pairs and the one-sided DtN map is not Hermitian, consistent with PML's need to absorb energy. On the other hand, our method propagates all the information from one end of the domain to the other end without any loss of energy. This is due to the fact that, even though our stiffness matrix is complex symmetric, the two-sided DtN map for the entire mesh is real symmetric (Hermitian) which is proved in Appendix A. Hermitian two-sided DtN map is partly a consequence of complex conjugate pairs of element lengths, which is consistent with keeping the boundary points invariant (which are on the real line).

2. Generalized eigenvalues of the stiffness matrix with respect to the mass matrix, in spite of both being complex-valued are non-negative and real. This property which is proved in Appendix B has implications in solving the eigenvalue problems. Also it is worthwhile noting that this phenomenon is linked to the important connection between Stieltjes rational approximants and relative rational approximants, recently found by Knizhnerman [18]. Specifically, he shows that there is a spectral equivalence between relative rational interpolants of the exponential and Stieltjes rational interpolants of the DtN map. As shown in [11], the latter corresponds to a spectrally matching three-point finite-difference scheme with symmetric negative operator, indicating that complex-length finite elements based on rational interpolants have real spectra. Rational interpolants include as a particular case the Padé approximant considered here, thus our finite-element scheme is exactly equivalent, in terms of the DtN maps, to the Padé finite-difference grid considered in [11]. Reference [11] also shows the equivalence of Padé finite-difference scheme and the spectral Galerkin method with polynomial basis functions, again in terms of the DtN maps.

While the discussion is limited to Laplace and Helmholtz equations, the mesh should also work for hyperbolic problems such as the wave equation. Since for the proposed discretization, the generalized eigenvalues of the stiffness matrix with respect to the mass matrix are non-negativeand real, traditional time-stepping methods would be stable. The only drawback of the proposed discretization for hyperbolic problems is that the mass matrix is not diagonal (due to midpoint integration). However, it could be block diagonal, as the degrees of freedom are coupled only in the $x$ direction; mass lumping can be performed in the $z$ direction. Fortunately, exponential convergence of the proposed method facilitates drastic reduction in the number of elements in the x direction, which translates into small block size and fairly efficient computation.

There may be concern associated with the increased cost due to complex arithmetic. This cost increase is negligible compared to the savings due to significantly reduced number of degrees of freedom facilitated by exponential convergence. Moreover, many problems involving wave propagation as well as electromagnetism require complex arithmetic and complex element lengths do not increase the computational cost for these problems.

The method in the presented form is applicable to problems where the solution is needed on the edges of the subdomains. If the solutions are needed only at the corner points, if the rest of the domain allows, the CFEM discretization can also be used in the transverse direction (if the material properties are pricewise constant in that direction), resulting in further reduction in the computational cost.

The current development revolves around scalar problems. Since fixed-point property of midpoint-integrated linear finite elements extends to vector equations (see [7]) it is expected that the proposed mesh would work equally well for electromagnetism and elasticity (similar to the development for finite difference schemes in [5, 6, 12]. We illustrate the applicability of the proposed approach for the special case of elastodynamics in the next section. Also successful application of the complex-length finite element method for solving the eigenvalue problems associated with elastic layered (stratified) media can be found in [19].

## 4. CFEM for vector equations (elastodynamics)

In this section we generalize the complex-length finite element formulation for solving linear elastodynamic problems. In particular, we prove the validity of the fixed point property for the vector elastic equations followed by driving the associated propagation factors. By observing that these propagation factors are identical to those in scalar analysis, we conclude that CFEM is applicable to vector wave equations.

*4.1. Model Problem*

We consider a two-dimensional elastic layered medium with in-plane deformation as shown in Fig. 1. Each layer is assumed to be homogeneous, but the material properties may vary between different layers. The equation representing in-plane wave propagation for the harmonic waves of the form $\bar{\mathbf{u}}(\mathbf{x},t) = \mathbf{u}(\mathbf{x},\omega)e^{-i\omega t} = \tilde{\mathbf{u}}(\mathbf{k},\omega)e^{i\mathbf{k}\cdot\mathbf{x}-i\omega t}$ with no external body forces and damping, can be written as,

$$\nabla_s^T \boldsymbol{\sigma} + \rho\omega^2 \mathbf{u} = \mathbf{0}, \tag{29}$$

where $\nabla_s^T = [\partial/\partial x \ 0 \ \partial/\partial z; \ 0 \ \partial/\partial z \ \partial/\partial x]$ is the gradient operator, $\omega \in \mathbb{R}$ denotes the temporal frequency, $\rho$ represents the density and $\mathbf{u} = \{u_x \ u_z\}^T$ is the infinitesimal in-plane displacement vector. Appropriate boundary conditions including Dirichlet, Neumann and Robin can be considered. The stress vector $\boldsymbol{\sigma} = \{\sigma_{xx} \ \sigma_{zz} \ \sigma_{xz}\}^T$ is related to the strain vector $\boldsymbol{\varepsilon} = \{\varepsilon_{xx} \ \varepsilon_{zz} \ \gamma_{xz}\}^T = \nabla_s \mathbf{u}$ by the stress–strain relationship:

$$\boldsymbol{\sigma} = \begin{bmatrix} D_{11} & D_{12} & D_{13} \\ D_{12} & D_{22} & D_{23} \\ D_{13} & D_{23} & D_{33} \end{bmatrix} \boldsymbol{\varepsilon}. \tag{30}$$

For the special case of isotropic elasticity, material coefficients are given in terms of Lamé constants $\lambda$ and $\mu$:

$$D_{11} = D_{22} = (\lambda + 2\mu), \ D_{33} = \mu, \ D_{12} = \lambda, \ D_{13} = D_{23} = 0. \tag{31}$$

Note that one can also model viscoelastic materials through frequency-dependent, complex Lamé constants. Expanding Eq. (29), we obtain,

$$\frac{\partial}{\partial x}\left(\mathbf{D}_{xx}\frac{\partial \mathbf{u}}{\partial x}\right)+\frac{\partial}{\partial x}\left(\mathbf{D}_{xz}\frac{\partial \mathbf{u}}{\partial z}\right)+\frac{\partial}{\partial z}\left(\mathbf{D}_{xz}^T\frac{\partial \mathbf{u}}{\partial z}\right)+\frac{\partial}{\partial z}\left(\mathbf{D}_{zz}\frac{\partial \mathbf{u}}{\partial z}\right)+\left(\rho\omega^2\mathbf{I}\right)\mathbf{u}=\mathbf{0}, \tag{32}$$

where $\mathbf{I}$ denotes $2\times 2$ identity matrix and $\mathbf{D}_{xx}$, $\mathbf{D}_{xz}$ and $\mathbf{D}_{zz}$ are the in-plane matrix coefficients that can be expressed in terms of the material property coefficients:

$$\mathbf{D}_{xx}=\begin{bmatrix} D_{11} & D_{13} \\ D_{13} & D_{33} \end{bmatrix},\quad \mathbf{D}_{xz}=\begin{bmatrix} D_{13} & D_{12} \\ D_{33} & D_{23} \end{bmatrix},\quad \mathbf{D}_{zz}=\begin{bmatrix} D_{33} & D_{23} \\ D_{23} & D_{22} \end{bmatrix}. \tag{33}$$

Discretizing Eq. (32) in the $z$ direction, we obtain,

$$\frac{\partial}{\partial x}\left(\mathbf{A}(z)\frac{\partial \mathbf{u}}{\partial x}\right)+\frac{\partial}{\partial x}\left(\mathbf{B}_1(z)\mathbf{u}\right)+\mathbf{B}_2(z)\frac{\partial \mathbf{u}}{\partial x}+\mathbf{D}(z,\omega)\mathbf{u}=\mathbf{0}, \tag{34}$$

where $\mathbf{A}$, $\mathbf{B}_1$, $\mathbf{B}_2$, and $\mathbf{D}$ are matrix differential operators of size $N\times N$, with $N$ being the number of degrees of freedom in the vertical direction. Eq. (34) admits individual modes of the form $\mathbf{u}=\boldsymbol{\phi}\,e^{ik_x x}$ where $k_x$ is the horizontal wavenumber. By substituting $\mathbf{u}$ in (34) we get the governing dispersion relation for $k_x$:

$$\left(-k_x^2\mathbf{A}+ik_x(\mathbf{B}_1+\mathbf{B}_2)+\mathbf{D}\right)\boldsymbol{\phi}=\mathbf{0}, \tag{35}$$

It can be shown that, considering an elastic layer from $x=0$ to $x=L$, the traction at $x=x_0$ takes the form,

$$\mathbf{F}\big|_{x=x_0}=-\left(\mathbf{A}\frac{\partial \mathbf{u}}{\partial x}+\mathbf{B}_1\mathbf{u}\right)\bigg|_{x=x_0}\equiv -\left(ik_x\mathbf{A}+\mathbf{B}_1\right)\mathbf{u}\big|_{x=x_0}, \tag{36}$$

Using the expression in (36), the displacement and stress on the left and right of the layer are given by:

$$\begin{aligned}\mathbf{u}_{x=0}^{\text{exact}}&=\mathbf{u}, & \mathbf{F}_{x=0}^{\text{exact}}&=-\left(ik_x\mathbf{A}+\mathbf{B}_1\right)\mathbf{u},\\ \mathbf{u}_{x=L}^{\text{exact}}&=e^{ik_x L}\mathbf{u}, & \mathbf{F}_{x=L}^{\text{exact}}&=-e^{ik_z L}\left(ik_x\mathbf{A}+\mathbf{B}_1\right)\mathbf{u}.\end{aligned} \tag{37}$$

*4.2. Fixed-point property and propagation factor of the mid-point integrated element*

We now focus on the segment $x=0$ to $x=x_1=L$. The fixed point property as defined in Section 2.3 is the recovery of the exact stiffness of half-space $[0,\infty)$ after augmenting the exact stiffness of the element $[0,L]$ with the exact stiffness of right half-space $[L,\infty)$. This definition is equivalent to satisfying the following equation, similar to (12):

$$\begin{bmatrix} \mathbf{K}_{11} & \mathbf{K}_{12} \\ \mathbf{K}_{21} & \mathbf{K}_{22}+\mathbf{K}_{HS} \end{bmatrix}\begin{Bmatrix} \mathbf{u}_0 \\ \mathbf{u}_L \end{Bmatrix}=\begin{Bmatrix} \mathbf{K}_{HS}\mathbf{u}_0 \\ \mathbf{0} \end{Bmatrix}, \tag{38}$$

where $\mathbf{u}_0$ and $\mathbf{u}_L$ are the displacements at $x=0$ and $x=L$, respectively. $\mathbf{K}_{kl}$ ($k,l \in \{1,2\}$) correspond to the stiffness components of the layer and,

$$\mathbf{K}_{HS} = -(ik_x \mathbf{A} + \mathbf{B}_1), \tag{39}$$

is the (exact) stiffness of the right half-space (see (36)).

Drawing from the observations from scalar problem, it is obvious that the fixed point property will not be satisfied if the layer stiffness $[0,L]$ is approximated by a regular finite element. However, like in the case of scalar problem, it turns out that the proposed mid-point integrated layers have the fixed point property and can recover the half-space stiffness exactly. This can be proved by verifying (38) for the discretized layer stiffness using midpoint integration with the following form,

$$\begin{bmatrix} \mathbf{K}_{11} & \mathbf{K}_{12} \\ \mathbf{K}_{21} & \mathbf{K}_{22} \end{bmatrix} = \frac{1}{L}\begin{bmatrix} \mathbf{A} & -\mathbf{A} \\ -\mathbf{A} & \mathbf{A} \end{bmatrix} + \frac{1}{2}\begin{bmatrix} -\mathbf{B}_1 & -\mathbf{B}_1 \\ \mathbf{B}_1 & \mathbf{B}_1 \end{bmatrix} - \frac{1}{2}\begin{bmatrix} -\mathbf{B}_2 & \mathbf{B}_2 \\ -\mathbf{B}_2 & \mathbf{B}_2 \end{bmatrix} - \frac{L}{4}\begin{bmatrix} \mathbf{D} & \mathbf{D} \\ \mathbf{D} & \mathbf{D} \end{bmatrix}. \tag{40}$$

Eq. (38) can be equivalently expressed as,

$$\begin{bmatrix} \mathbf{K}_{11} - \mathbf{K}_{HS} & \mathbf{K}_{12} \\ \mathbf{K}_{21} & \mathbf{K}_{22} + \mathbf{K}_{HS} \end{bmatrix} \begin{Bmatrix} \mathbf{u}_0 \\ \mathbf{u}_L \end{Bmatrix} = \begin{Bmatrix} \mathbf{0} \\ \mathbf{0} \end{Bmatrix}. \tag{41}$$

Utilizing (39) and (40) in (41) yields:

$$\left( \frac{1}{L}\begin{bmatrix} (1+ik_x L)\mathbf{A} & -\mathbf{A} \\ -\mathbf{A} & (1-ik_x L)\mathbf{A} \end{bmatrix} + \frac{1}{2}\begin{bmatrix} (\mathbf{B}_1+\mathbf{B}_2) & -(\mathbf{B}_1+\mathbf{B}_2) \\ (\mathbf{B}_1+\mathbf{B}_2) & -(\mathbf{B}_1+\mathbf{B}_2) \end{bmatrix} - \frac{L}{4}\begin{bmatrix} \mathbf{D} & \mathbf{D} \\ \mathbf{D} & \mathbf{D} \end{bmatrix} \right) \begin{Bmatrix} \mathbf{u}_0 \\ \mathbf{u}_L \end{Bmatrix} = \begin{Bmatrix} \mathbf{0} \\ \mathbf{0} \end{Bmatrix}. \tag{42}$$

As $\mathbf{u}_0$ and $\mathbf{u}_L$ satisfy the dispersion relation (35) we have,

$$(\mathbf{B}_1 + \mathbf{B}_2)\mathbf{u}_0 = -(ik_x \mathbf{A} + \mathbf{D}/ik_x)\mathbf{u}_0 \quad \text{and} \quad (\mathbf{B}_1 + \mathbf{B}_2)\mathbf{u}_L = -(ik_x \mathbf{A} + \mathbf{D}/ik_x)\mathbf{u}_L. \tag{43}$$

Substituting these into Eq. (42) and pre-multiplying by the nonsingular matrix $[L/2 \ -L/2; \ -ik_x \ -ik_x]$ gives:

$$\begin{bmatrix} (1+ik_x L/2)\mathbf{A} & (-1+ik_x L/2)\mathbf{A} \\ (1+ik_x L/2)\mathbf{D} & (-1+ik_x L/2)\mathbf{D} \end{bmatrix} \begin{Bmatrix} \mathbf{u}_0 \\ \mathbf{u}_L \end{Bmatrix} = \begin{Bmatrix} \mathbf{0} \\ \mathbf{0} \end{Bmatrix}, \tag{44}$$

which is satisfied with,

$$\mathbf{u}_L = P\,\mathbf{u}_0 = \left( \frac{1+ik_x L/2}{1-ik_x L/2} \right) \mathbf{u}_0, \tag{45}$$

where $P$ is the propagation factor of the displacement using one element. So far we showed the recovery of the exact half-space stiffness after adding a mid-point integrated linear finite element to a half-space, *irrespective* of the element length $L$. Besides we obtained the propagation factor of the displacement field for the mid-point integrated element as given in (45). While this derivation for displacement propagation factor is also presented in an earlier paper [7], CFEM validity requires the propagation factor for displacement-traction pairs similar to the scalar case. We now proceed to derive the propagation factor of the traction (or equivalently force).

Given the displacement and stress at the beginning of the segment, we need to find the CFEM approximation of the displacement and stress at the end. To this end we look into the relation between the traction and displacement for the mid-point integrated element,

$$\begin{Bmatrix} \mathbf{F}_0 \\ -\mathbf{F}_L \end{Bmatrix} = \begin{bmatrix} \mathbf{K}_{11} & \mathbf{K}_{12} \\ \mathbf{K}_{21} & \mathbf{K}_{22} \end{bmatrix} \begin{Bmatrix} \mathbf{u}_0 \\ \mathbf{u}_L \end{Bmatrix}, \tag{46}$$

where the stiffness matrix of the element is given in (40). Using the fixed point property given in (41), the above expression can be written in the form,

$$\begin{Bmatrix} \mathbf{F}_0 \\ -\mathbf{F}_L \end{Bmatrix} = \left( \begin{bmatrix} \mathbf{K}_{11} - \mathbf{K}_{HS} & \mathbf{K}_{12} \\ \mathbf{K}_{21} & \mathbf{K}_{22} + \mathbf{K}_{HS} \end{bmatrix} + \begin{bmatrix} \mathbf{K}_{HS} & \mathbf{0} \\ \mathbf{0} & -\mathbf{K}_{HS} \end{bmatrix} \right) \begin{Bmatrix} \mathbf{u}_0 \\ \mathbf{u}_L \end{Bmatrix} \stackrel{\text{Fixed-Point}}{=} \begin{bmatrix} \mathbf{K}_{HS} & \mathbf{0} \\ \mathbf{0} & -\mathbf{K}_{HS} \end{bmatrix} \begin{Bmatrix} \mathbf{u}_0 \\ \mathbf{u}_L \end{Bmatrix}, \tag{47}$$

which gives the traction $\mathbf{F}_0$ and $\mathbf{F}_L$ on the left and right side of the elements as,

$$\mathbf{F}_0 = -\left(ik_x \mathbf{A} + \mathbf{B}_1\right)\mathbf{u}_0, \quad \mathbf{F}_L = -P\left(ik_x \mathbf{A} + \mathbf{B}_1\right)\mathbf{u}_0 = P\mathbf{F}_0. \tag{48}$$

Considering (45) and (48) indicates that the propagation factors are the same for displacements and tractions. Further, comparing these with (37) indicates that the exponential propagator in (37) is approximated by the propagator $P$. Utilizing (45), we have,

$$e^{ik_x L} \approx P = \left( \frac{1 + ik_x L/2}{1 - ik_x L/2} \right). \tag{49}$$

The above expression is identical to the propagator in (25) which associated with the scalar equation. Given this identical behavior of midpoint integrated linear finite elements between scalar and vector equations, the CFEM construction for scalar equations in after (25) is immediately applicable for the vector wave equation.

## 5. Numerical examples

### 5.1. Two-point boundary value problem

In this example, we consider a simple two-point boundary value problem on a unit interval with homogeneous Dirichlet condition applied on the right. Unit Neumann data is specified on the left boundary and the goal is to obtain the displacements at the boundary. The relative error is computed using the following norm:

$$error = \frac{|u_{exact} - u_{approx}|}{|u_{exact}| + |u_{approx}|}. \tag{50}$$

The analysis is performed for elliptic version of one-dimensional Eq. (8) with positive $\lambda = k^2$ with $10 \leq k \leq 200$. CFEM meshes with increasing number of elements (ranging from 1 through 40) are used to test the convergence of the method. Fig. 3 shows the results from the analysis, where various curves correspond to various values of $k$. It can be seen that the error converges until it reaches the roundoff limit of approximately $10^{-16}$. Furthermore, the convergence is super-exponential, similar to the convergence of the optimal grids proposed in [4]. More quantitatively, even for a large value of $k = 200$, where the solution contains steeply varying functions $e^{\pm 200}$, a coarse mesh of just 20 elements can get the relative error down to $10^{-4}$, clearly illustrating the efficiency of the proposed discretization.

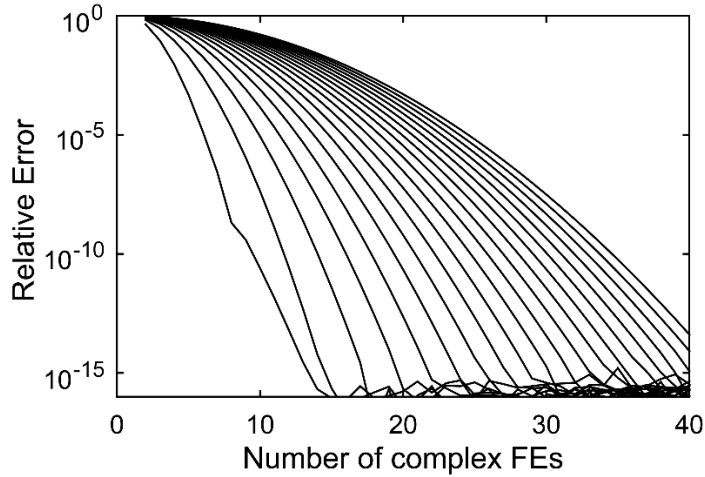

**Fig. 3.** Convergence for the two-point elliptic boundary value problem. The left most curve corresponds to $k = 10$, while the right most curve corresponds to an extremely high decay parameter of $k = 200$.

*5.2. Two-point boundary value problem: Helmholtz equation*

In order to assess the performance of the proposed method for wave propagation problems, we turn to the Helmholtz equation, which has a negative value of $\lambda = -\omega^2$, with the frequency range of $4 \leq \omega \leq 40$. As in the previous numerical experiment, homogeneous Dirichlet condition is applied on the right and unit Neumann data is specified on the left. The error as defined in (50) is plotted in Fig. 4. Super-exponential convergence is again observed, but the convergence appears to occur only after the number of elements exceeds a threshold. This threshold corresponds to an average of three elements per wavelength. This is expected because, our finite element discretization is equivalent at the boundary points to Spectral Legendre-Galerkin method (see Section 3.4), and the latter has similar threshold; a more efficient approximant, e.g. Pade-Chebyshev, would make this threshold very close to the Nyquist limit of two points per wavelength [3, 4, 6]. A peculiar behavior observed in Fig. 4 is that, as the frequency is increased, the error does not converge to zero, but to a nonzero value that grows as the frequency is further increased. It turns out that reordering of elements alleviates this problem, as described below.

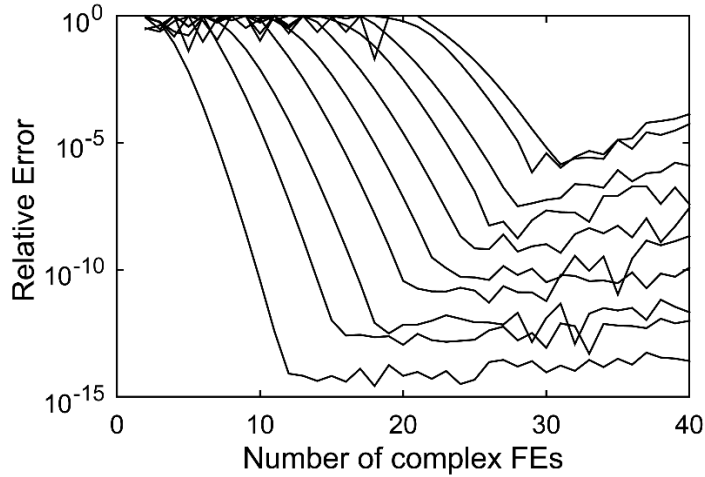

**Fig. 4.** Convergence for the two-point Helmholtz boundary value problem for varying frequencies. The left most curve in each figure corresponds to the lowest frequency of 4.0, while the right most curve corresponds to the highest frequency of 40. The meshes shown in Fig. 2 are used for computation. Note that the error does not converge to a small value, especially for high frequencies.

*Element reordering*: An intuitive argument for the numerical problems is as follows. Consider a wave of the form $u = ae^{\pm ikx}$ propagating through the meshes shown in Fig. 2. As we approach the center of the mesh, the imaginary part of $x$ increases, indicating that $u$ could grow exponentially. While $u$ would eventually decay back to its original order of magnitude as we approach the right edge, the growth at the center could be very large, especially for high frequencies. We hypothesize that this translates into growth of the round-off error, resulting in significant numerical inaccuracies. An option to alleviate this problem is to pair the elements with complex conjugate lengths. Such pairing minimizes the artificial growth of the field variable and thus the round-off error. A disadvantage of this approach is that the symmetry in the domain geometry is not preserved after discretization. To achieve mesh symmetry, we propose to swap every alternate element in the left half of the mesh in Fig. 2, with its symmetric counterparts on the right (few examples of rearranged meshes are given in Fig. 5). This enforces alternating signs for imaginary parts of element lengths, thus reducing the growth of the imaginary part of the coordinate $x$.

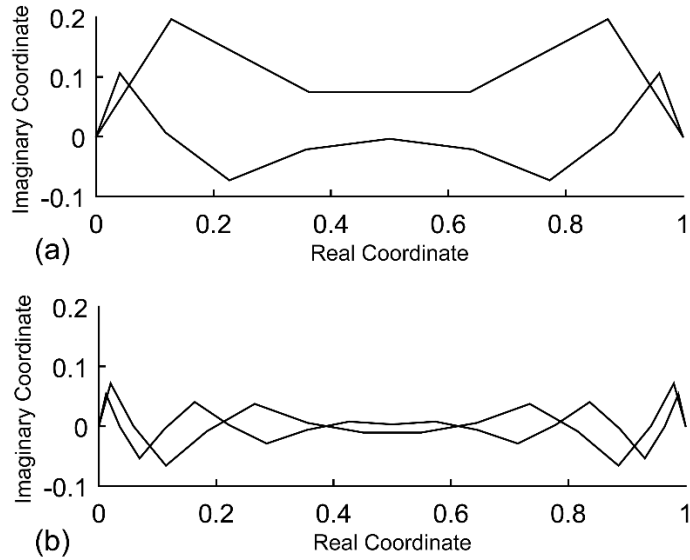

**Fig. 5.** Representative meshes after element reordering for (a) 5-element and 10-element meshes and (b) 15-element and 20-element meshes. Note that the bounds of imaginary part of the nodal coordinates are reducing with refinement. This counters the numerical growth and helps achieve better convergence (compare these meshes with the meshes in Fig. 2).

With the new element ordering, the solution to the Helmholtz equation improved considerably, as shown in Fig. 6. Robust mathematical analysis of the effects of ordering on numerical accuracy would be considered in the future.

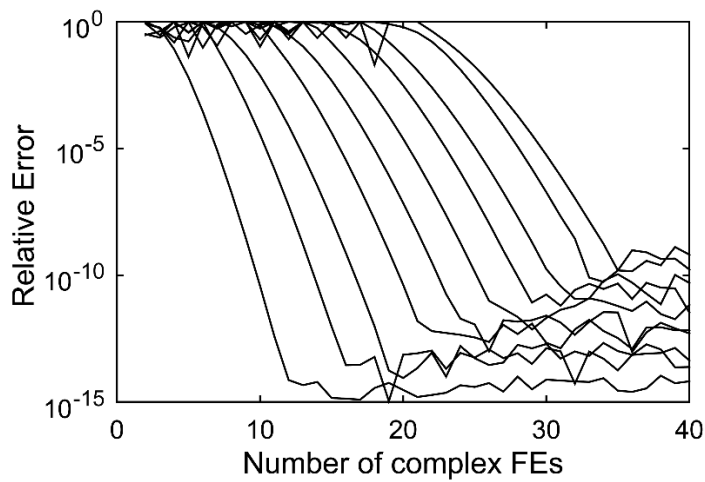

**Fig. 6.** Convergence for the two-point Helmholtz boundary value problem after element reordering shown in Fig. 5. In comparison with Fig. 4, the error is converging to a much smaller value, indicating the effectiveness of element reordering.

In order to investigate the efficiency of the CFEM, we also compare the convergence results of CFEM with spectral finite elements (SFEM) using Lobatto polynomials as discussed in [20]. Results are reported in Fig. 7 for the lowest and highest values of $k$ and $\omega$ for the elliptic and Helmholtz two-point boundary value problems, respectively.

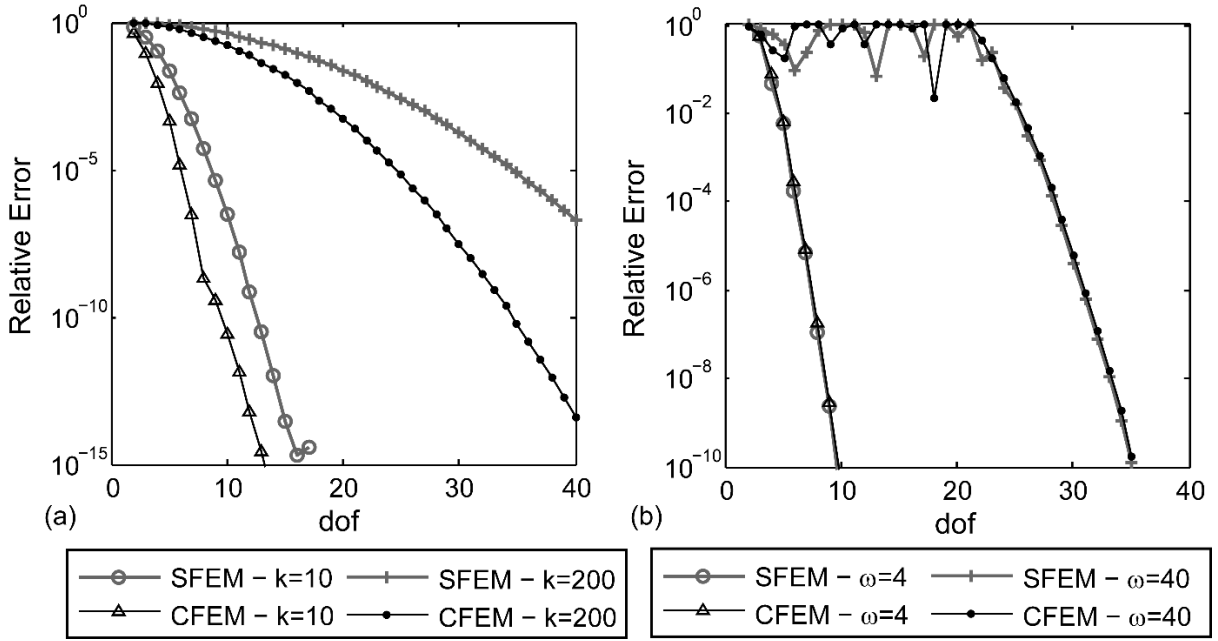

**Fig. 7.** Convergence curves for (a) the two-point elliptic boundary value problem and (b) the two-point Helmholtz boundary value problem using complex-length and spectral finite elements

Fig. 7 (a) shows that for the elliptic equation, performance of CFEM is slightly better than SFEM for $k = 10$ and the disparity between the two methods becomes larger with increasing $k$. However as depicted in Fig. 7 (b) for the Helmholtz equation which has a harmonic response, both methods have the same convergence for different frequencies. It is important to recall that even in the case that the number of total degrees of freedom is the same by adopting either SFEM or CFEM for reaching to a given accuracy, CFEM matrices are sparse (tridiagonal in this case), whereas SFEM leads to full matrices. CFEM would thus have reduced computational cost by adopting efficient sparse solution techniques, e.g. block tri-diagonal solvers in [21]. Furthermore, CFEM can easily be implemented in existing finite element codes with simple change of element lengths and midpoint integration.

*5.3. Two-dimensional layer: Laplace equation*

Extending the one-dimensional examples to two-dimensional setting, we consider a long layer of size

$1\times10$, governed by Laplace equation with unit coefficients. Homogeneous Dirichlet boundary condition is applied at the bottom and homogeneous Neumann condition is applied on the top and right edges. In order to focus on the effect of $x$ discretization, we use a fine mesh of 200 elements in the $z$ direction. The excitation (specified Neumann data) at the left boundary is given by the analytic function

$$-G\frac{\partial u}{\partial x} = e^{16+4/y/(y-1)} \,. \tag{51}$$

Such excitation is chosen to excite various vertical harmonics and thus various decay rates in the horizontal direction. In order to compute the error associated with various analyses, the converged FEM solution is considered here as the reference solution. The response on the left edge is computed and the resulting error $e_r = \|\mathbf{u}-\mathbf{u}^{ref}\|_2 / \|\mathbf{u}^{ref}\|_2$ is examined, where $\mathbf{u}$ and $\mathbf{u}^{ref}$ are the calculated and reference displacement vectors, respectively.

The convergence of the error with mesh refinement using CFEM is shown in Fig. 8. Also for the sake of comparison, convergence analysis is performed for regular finite element method and the results are shown in Fig. 8 (b). Clearly CFEM shows exponential convergence, while regular finite elements show expected second-order convergence (note that the left plot is in semi-log scale, while the right plot is in log-log scale). Quantitatively, just 10 complex-length finite elements are sufficient to achieve an error of less than 1%, whereas uniform finite element mesh requires more than 100 elements for the same accuracy. When the desired error is reduced 0.01%, we require just 14 complex-length finite elements, as opposed to approximately 1000 regular finite elements.

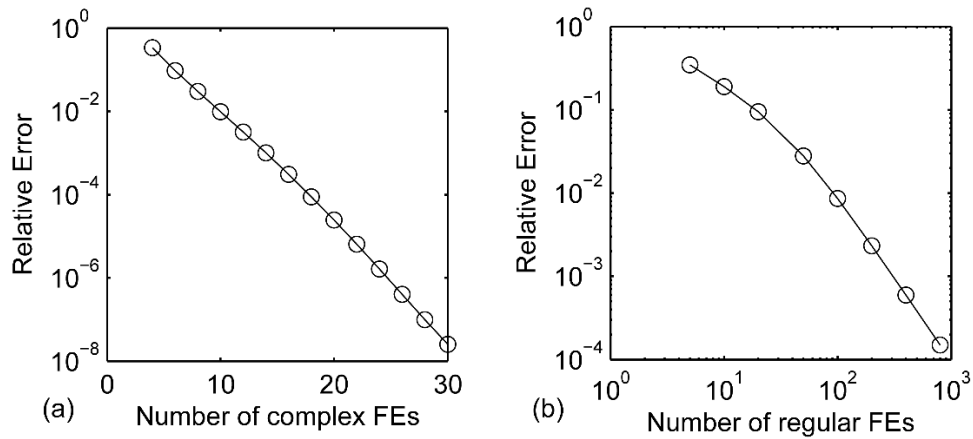

**Fig. 8.** Convergence of complex-length finite element discretization for Laplace equation in two-dimensional domain. Note the exponential convergence in complex-length finite element discretization as opposed to algebraic convergence of regular FEM. Also note that, for 1% relative error, one would need 10 complex-length finite elements as opposed to 100 regular elements. 0.01% relative error requires just 14 complex-length finite elements, as opposed to 1000 regular finite elements.

*5.4. Two-dimensional layer: Helmholtz equation*

The 2-D problem of the previous section is solved for Helmholtz equation with frequency of $\omega=3$. The frequency is chosen to ensure existence of significant evanescent as well as propagating waves. Since the domain is bounded, in order to eliminate any possible complications due to resonance, we used $G=1+0.01i$ in the governing equation (1), which has the effect of damping out the resonance. Convergence analysis similar to that performed for the Laplace equation is performed for this problem. The results are shown in Fig. 9 (a). The plot also contains the convergence results for the field variable at $x=10$. Possibly because of resonances and oscillations, the convergence is slowed compared to the Laplace equation. However, it is much faster compared to regular finite element discretization, as described in the following paragraph.

The convergence analysis is also performed for regular finite element discretization and the results are shown in Fig. 9 (b). As expected, the convergence is second-order. Quantitatively, a desired error of 1% requires just 17 complex-length finite elements, as opposed to 400 regular FE elements. When the desired error is 0.1%, the disparity increases: 20 elements for complex-length finite element discretization as opposed to more than 1000 elements for regular FE discretization (even when 1000 elements are used in the $x$ direction, the element size is $0.1\times 0.05$, which explains why the error does not reach the asymptote of the $z$ discretization error). Note that, due to exponential convergence, complex-length finite element discretization does not suffer from adverse consequences of pollution effects of the dispersion error associated with regular finite element discretization [22].

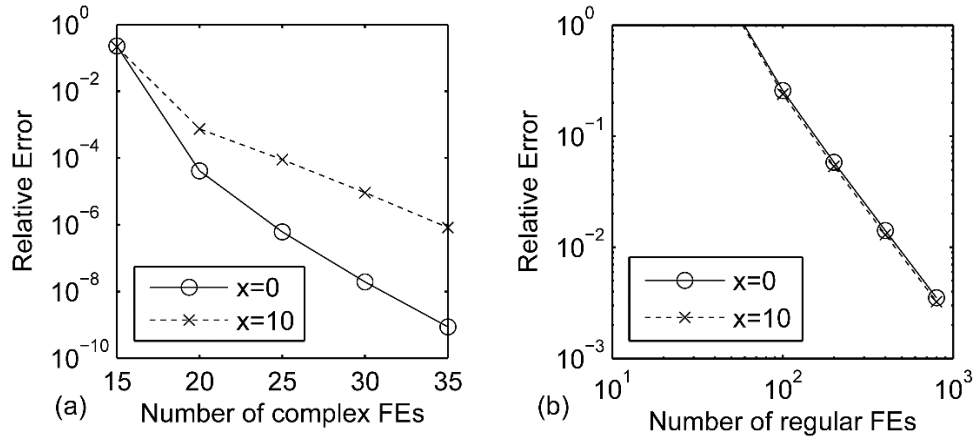

**Fig. 9.** Convergence of complex-length finite element discretization for Helmholtz equation in a two-dimensional setting. Note the exponential convergence in complex FEM. Also, note that for 1% relative error, we would need 17 complex-length finite elements as opposed to 400 regular finite elements. 0.1% relative error requires just 20 complex-length finite elements, as opposed to more than 1000 regular finite elements.

*5.5. Two-dimensional multi-domain layer: Helmholtz equation*

To illustrate the applicability of the proposed method for multi-domain problems, the modulus of the right half $(5 < x < 10)$ is modified as $G = 2 + 0.02i$. The sub-domains $(0,5)$ and $(5,10)$ are discretized separately with complex-length finite elements, with nodes shared at $x = 5$. Note that accuracy can be expected at $x = 0, 5, 10$, while the remaining nodes have complex coordinates. The convergence of the response at these locations is plotted in Fig. 10 (a). Similar to the previous example, the convergence is worse compared to the Laplace equation, but is much faster than that for regular finite element discretization, which is shown in Fig. 10 (b).

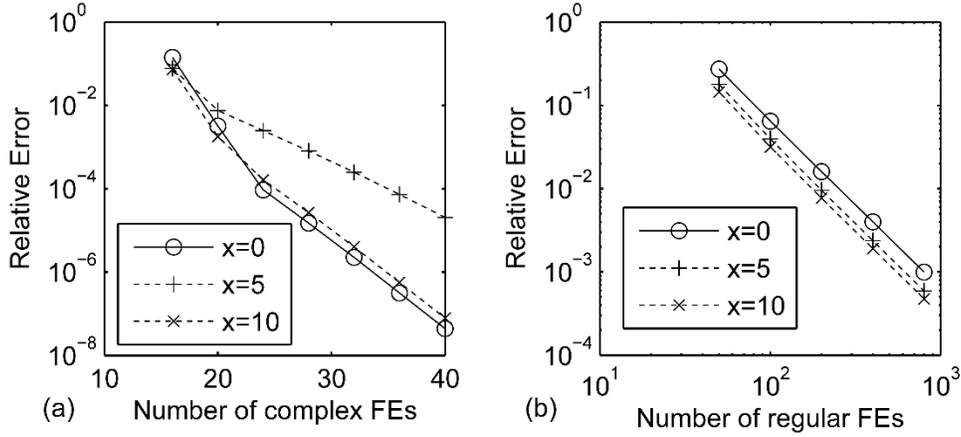

**Fig. 10.** Convergence for multidomain problem. Note the exponential convergence in complex FEM. For 1% relative error, we would need 20 complex-length finite elements as opposed to 300 regular finite elements. 0.1% error tolerance requires 28 complex-length finite elements as opposed to more than 800 regular finite elements

*5.6. Two-dimensional multi-domain layer: time-harmonic elastodynamics*

The 2-D multi-domain problem of the previous section is now solved for time-harmonic elastodynamic equation given in (29) with frequency $\omega = 3$. The shear modulus of the left half and right half are chosen as $G = 1 + 0.01i$ and $G = 2 + 0.02i$, respectively. The Poisson's ratio is taken as $\nu = 0.35$, the density is assumed to be $\rho = 1$ and the following horizontal traction is applied on the left boundary of the domain at $x = 0$:

$$\left\{ \begin{matrix} \sigma_{xx} \\ \sigma_{xy} \end{matrix} \right\}\bigg|_{x=0} = -\left( \mathbf{G}_{xx}\frac{\partial}{\partial x} + \mathbf{G}_{xy}\frac{\partial}{\partial y} \right)\mathbf{u}\bigg|_{x=0} = \left\{ \begin{matrix} e^{16 + 4/y(y-1)} \\ 0 \end{matrix} \right\}. \tag{52}$$

Similar to the previous example the displacement convergence along the interfaces $x = 5, 10, 15$ are shown in Fig. 11 using complex and regular finite elements with the relative error including both horizontal and vertical components of the displacement. It can be observed that, like in the previous example, using CFEM requires much fewer elements when compared to regular finite elements.

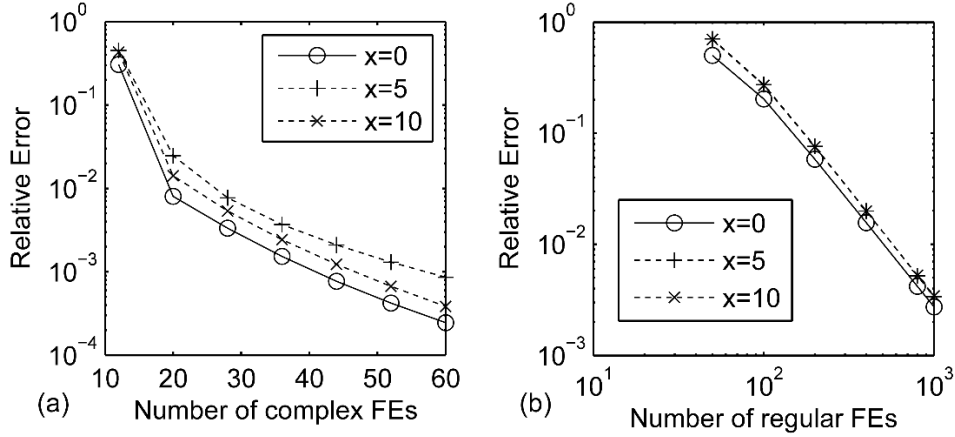

**Fig. 11.** Convergence for multidomain elastic problem. For 1% relative error, we would need 28 complex-length finite elements as opposed to 500 regular finite elements. 0.1% error tolerance requires 60 complex-length finite elements as opposed to more than 1000 regular finite elements.

## 6. Conclusions

We introduced a novel modification of linear finite element discretization to achieve exponential convergence of the solution at select points in the domain. The main idea is to obtain exponentially convergent approximation for the DtN map of the sub-domains spanning between the points of interest. By employing midpoint integration rules and an unconventional bending of the finite element mesh into the complex space, we are able to achieve high accuracy at the edges of the sub-domain with a very coarse discretization of the interior. The development is facilitated by linking midpoint-integrated linear finite element discretization to Crank-Nicolson stepping of the associated first order form, which in turn is related to Padé approximants of the exponential function. The parameters of the Padé approximants are translated back to the parameters of finite element discretization (element lengths). The resulting element lengths are complex-valued, and the method is named complex-length finite element method. The method inherits the exponential convergence from the underlying Padé approximants. Exponential convergence is verified with a variety of numerical experiments involving Laplace, Helmholtz and elastodynamic equations. The examples also indicate that, for practical error tolerance requirements, the method facilitates reduction of the number of elements by an order of magnitude resulting in significant savings in computational cost.

The current paper indicates promise of the complex-length finite element method, but further investigations and enhancements are necessary, especially because of the unconventional nature of the method. These investigations include: (a) deeper mathematical understanding of the underlying approximation properties; (b) investigation of numerical properties of the resulting discretizations; (c) extension to tensor product complex-valued discretizations for higher dimensions and resolution of possible stability issues associated with midpoint integration in higher dimensions; (d) extension of the method for hyperbolic problems; and (e) possible extension to topologically triangular meshes thus to unstructured meshes. These issues are subjects of ongoing investigation.

**Acknowledgements**

The authors wish to thank Leonid Knizhnerman for useful discussions. The work is partially supported by National Science Foundation (Grant No. DMS-1016514) and Schlumberger Technology Corporation.

**Appendix A: Proof of DtN map being Hermitian**

We consider the following Equation Similar to Eq. (8) with $\lambda = k^2$,

$$-\frac{\partial^2 u}{\partial x^2} + k^2 u = 0, \quad \text{for} \quad 0 < x < L. \tag{A1}$$

Discretizing the weak form as presented in (15), dynamic stiffness matrix of any mid-point integrated element $(\mathcal{K}_{\text{elem}} + k^2 \mathcal{M}_{\text{elem}})$, similar to Eq. (17) with $\lambda = k^2$, can be written in the form $\mathcal{K}_j^d = [a, \ b; \ b \ a]$ where $a = 1/L_j + k^2 L_j / 4$ and $b = -1/L_j + k^2 L_j / 4$. Therefore, for a given $k^2 \in \Re$, the assembled stiffness matrix for a pair of elements with complex conjugate lengths $L$ and $\bar{L}$ is,

$$(\mathcal{K}_{L,\bar{L}} + k^2 \mathcal{M}_{L,\bar{L}}) = \begin{bmatrix} a & b & \\ b & a+\bar{a} & \bar{b} \\ & \bar{b} & \bar{a} \end{bmatrix}. \tag{A2}$$

By condensing out the middle node and using the fixed point property, i.e., $a^2 - b^2 = k^2$ (see (14)), the $2 \times 2$ DtN map of (A2) takes the form,

$$\mathbf{S}_{L,\bar{L}} = \frac{1}{2\Re(a)} \begin{bmatrix} k^2 + |a|^2 & -|b|^2 \\ -|b|^2 & k^2 + |a|^2 \end{bmatrix}. \tag{A3}$$

This indicates that $\mathbf{S}_{L,\bar{L}}$ is a real symmetric (Hermitian) matrix similar to the stiffness matrix of an element with real length. Since the CFEM mesh consists of either pairs of elements with complex conjugate lengths or an element with real length, the assembly of DtN maps for all elements and consequently the DtN map of the entire mesh will be real symmetric (Hermitian). This is true not only for positive $k^2$ (elliptic equation) but also for negative $k^2$ (Helmholtz equation).

**Appendix B: Proof of eigenvalues of $(\mathcal{K}, \mathcal{M})$ being real and nonnegative**

Discretizing the weak form of Eq. (A1) leads to the system of equations $(\mathcal{K}+k^2\mathcal{M})$. Generalized eigenvalues of the pencil $(\mathcal{K}, \mathcal{M})$ are real and non-negative, if and only if all values of $k$ with singular $(\mathcal{K}+k^2\mathcal{M})$ are purely imaginary. This will be proved using the following lemmas.

**Lemma 1.** For a mesh with elements $L_1,...,L_n$, the propagator $P_{L_1,...,L_n}(k) = \pm 1$ for any $k$ with singular $(\mathcal{K}+k^2\mathcal{M})$.

**Proof.** Eq. (17) shows that any mid-point integrated linear element in the mesh satisfies the fixed-point property. Therefore all the elements $L_1,...,L_n$ and consequently the two-point DtN map of the entire mesh will satisfy this property. Accordingly the two-point DtN of the entire mesh will take the form,

$$\mathbf{S} = \begin{bmatrix} S_{\text{diag}} & S_{\text{off}} \\ S_{\text{off}} & S_{\text{diag}} \end{bmatrix}. \tag{B1}$$

Note that the diagonal entries of DtN map are the same due to the mirror symmetry of the mesh (note that $\mathbf{S}$ is invariant of element ordering as discussed in Section 3.3 and we simply choose the mesh that is mirror symmetric to make the claim). Similar to (12), the fixed-point property of this DtN map can be written as follows,

$$\begin{bmatrix} S_{\text{diag}} - K_{\text{halfspace}} & S_{\text{off}} \\ S_{\text{off}} & S_{\text{diag}} + K_{\text{halfspace}} \end{bmatrix} \begin{Bmatrix} u_0 \\ u_1 \end{Bmatrix} = \begin{Bmatrix} 0 \\ 0 \end{Bmatrix}. \tag{B2}$$

Given that $K_{\text{halfspace}} = k$ and using the propagator for $u_1 = P_{L_1,...,L_n} u_0$, (B2) becomes,

$$\begin{bmatrix} S_{\text{diag}} - k & S_{\text{off}} \\ S_{\text{off}} & S_{\text{diag}} + k \end{bmatrix} \begin{Bmatrix} 1 \\ P_{L_1,...,L_n} \end{Bmatrix} = \begin{Bmatrix} 0 \\ 0 \end{Bmatrix}. \tag{B3}$$

Using (B3) and the fixed point property $S_{\text{diag}}^2 - S_{\text{off}}^2 = k^2$ (see Eq. (14)), the propagator $P_{L_1,...,L_n}$ can be obtained as,

$$P_{L_1,...,L_n} = \pm\sqrt{\frac{1 - k/S_{\text{diag}}}{1 + k/S_{\text{diag}}}}. \tag{B4}$$

Note that the singularity of $(\mathcal{K}+k^2\mathcal{M})$ implies the singularity (rank-deficiency) of $\mathbf{S}$. However due to fixed-point property $\det(\mathbf{S}) = S_{\text{diag}}^2 - S_{\text{off}}^2 = k^2$, which is finite. This can only occur if the components of the DtN map, $S_{\text{diag}}, S_{\text{off}} \to \infty$. With this, (B4) becomes,

$$P_{L_1,\ldots,L_n}(k) = \pm 1. \tag{B5}$$

An alternative proof of Lemma 1 can be obtained by utilizing the isomorphism between the propagators obtained via CFEM and the two-sided optimal grid in Section 3.3 of [23], but is not presented here.

**Lemma 2.** For a complex $k = k_R + ik_I$ with nonzero $k_R$, the absolute value of the propagator cannot be unity $|P_{L_1,\ldots,L_n}(k)| \neq 1$.

**Proof.** As shown in Eq. (26) the propagator for the entire mesh with elements $L_1,\ldots,L_n$ can be written as the multiplication of the elements' individual propagators, i.e.,

$$P_{L_1,\ldots,L_n}(k) = P_{L_n} P_{L_{n-1}} \cdots P_{L_2} P_{L_1}. \tag{B6}$$

Figure 2 shows that CFEM mesh consists of either pairs of elements with complex conjugate lengths or an element with real length. We start with the propagator for a pair of elements with complex conjugate lengths $L$ and $\bar{L}$ (obtained using Eq. (26)):

$$P_{L,\bar{L}}(k) = \left(\frac{1+kL/2}{1-kL/2}\right)\left(\frac{1+k\bar{L}/2}{1-k\bar{L}/2}\right). \tag{B7}$$

Considering $L = L_R + iL_I$, the propagator can be written as,

$$P_{L,\bar{L}}(k) = \frac{(\mathcal{X}+L_R k_R)+i(\mathcal{Y}+L_R k_I)}{(\mathcal{X}-L_R k_R)+i(\mathcal{Y}-L_R k_I)}, \tag{B8}$$

where $\mathcal{X} = (k_R^2 - k_I^2)|L|^2/4 + 1$ and $\mathcal{Y} = k_R k_I |L|^2/2$. The absolute value takes the form,

$$|P_{L,\bar{L}}(k)| = \sqrt{\frac{\mathcal{Q}+\mathcal{E}}{\mathcal{Q}}}, \tag{B9}$$

where $Q = (\mathcal{X} - L_R k_R)^2 + (\mathcal{Y} - L_R k_I)^2$ and $\mathcal{E} = 4L_R(\mathcal{X}k_R + \mathcal{Y}k_I) = k_R L_R(|L|^2|k|^2 + 4)$. Now it can be seen that $k_R > 0$ makes $|P_{L,\bar{L}}(k)| > 1$ and $k_R < 0$ makes $|P_{L,\bar{L}}(k)| < 1$.

For odd number of complex finite elements, in addition to the complex conjugate lengths, there is a single element with real length. For this element with real length $L_R$, the associated propagator can be obtained by using (B7) with $L_I = 0$. This gives $P_{L_R}^2(k) = P_{L,\bar{L}}(k)$ and consequently $|P_{L,\bar{L}}(k)| = |P_{L_R}(k)|^2$. It can be easily seen that similar to complex conjugate element pairs, $k_R > 0$ makes $|P_{L_R}(k)| > 1$ and $k_R < 0$ makes $|P_{L_R}(k)| < 1$.

Now considering the entire propagator given in (B6) for a general mesh with $m$ pairs of elements with complex conjugate lengths and an element with real length,

$$\begin{cases} |P_{L_i,\bar{L}_i}|, |P_{L_R}| > 1 \quad (i=1,...,m) \Rightarrow |P_{L_1,...,L_n}| = |P_{L_R}| \times \prod_{i=1}^{m} |P_{L_i,\bar{L}_i}| > 1, & \text{for } k_R > 0, \\ |P_{L_i,\bar{L}_i}|, |P_{L_R}| < 1 \quad (i=1,...,m) \Rightarrow |P_{L_1,...,L_n}| = |P_{L_R}| \times \prod_{i=1}^{m} |P_{L_i,\bar{L}_i}| < 1, & \text{for } k_R < 0. \end{cases} \quad (B10)$$

From (B10) it can be concluded that $|P_{L_1,...,L_n}(k)| \neq \pm 1$

**Proof of the main claim:** Lemma 1 states that singularity of singular $(\mathcal{K} + k^2 \mathcal{M})$ implies $P_{L_1,...,L_n}(k) = \pm 1$. On the other hand, from Lemma 2, any $k$ that is not purely imaginary would result in $P_{L_1,...,L_n}(k) \neq \pm 1$, which indicates $P_{L_1,...,L_n}(k) = \pm 1$ is possible for purely imaginary $k$. Thus, singularity of $(\mathcal{K} + k^2 \mathcal{M})$ implies that $k$ must be purely imaginary, which in turn implies that the generalized eigenvalues of $\mathcal{K}$ with respect to $\mathcal{M}$ are real and nonnegative.